\documentclass{amsart}

\usepackage{amssymb}

\newtheorem{thm}{Theorem}

\newtheorem{lem}[thm]{Lemma}
\newtheorem{cor}[thm]{Corollary}

\newtheorem{prop}[thm]{Proposition}

\theoremstyle{definition}
\newtheorem{defn}[thm]{Definition}

\newtheorem{say}[thm]{}
\newtheorem{exmp}[thm]{Example}

\newtheorem{rem}[thm]{Remark}          
\newtheorem{note}[thm]{Note}            
\newtheorem{ack}{Acknowledgments}        
\newtheorem{notation}[thm]{Notation}   
  
\newtheorem{defn-thm}[thm]{Definition--Theorem}  
\newtheorem{defn-lem}[thm]{Definition--Lemma}  

\theoremstyle{remark}
\newtheorem{claim}[thm]{Claim}

\renewcommand{\c}[0]{{\mathbb C}}  

\renewcommand{\o}[0]{{\mathcal O}} 
\newcommand{\z}[0]{{\mathbb Z}}

\renewcommand{\a}[0]{{\mathbb A}}

\newcommand{\p}[0]{{\mathbb P}}

\newcommand{\q}[0]{{\mathbb Q}}

\newcommand{\qtq}[1]{\quad\mbox{#1}\quad}
\newcommand{\spec}[0]{\operatorname{Spec}}
\newcommand{\pic}[0]{\operatorname{Pic}}

\newcommand{\supp}[0]{\operatorname{Supp}}    
\newcommand{\red}[0]{\operatorname{red}}

\newcommand{\aut}[0]{\operatorname{Aut}}

\newcommand{\sing}[0]{\operatorname{Sing}}

\newcommand{\chr}[0]{\operatorname{char}}

\newcommand{\weil}[0]{\operatorname{Weil}}

\newcommand{\cl}[0]{\operatorname{Cl}}

\newcommand{\rdown}[1]{\lfloor{#1}\rfloor}

\newcommand{\onto}[0]{\twoheadrightarrow}

\newcommand{\lcm}[0]{\operatorname{lcm}}

\def\into{\DOTSB\lhook\joinrel\to}

\begin{document}
\bibliographystyle{amsalpha}

\title{Seifert $G_m$-bundles}
\author{J\'anos Koll\'ar}
\maketitle

Seifert fibered 3--manifolds were introduced and studied in \cite{seif}.
(See  \cite{seif-book} for an English translation.)
Roughly speaking, these are 3--manifolds $M$ which admit a
differentiable map $f:M\to F$ to a surface $F$
such that every fiber is a circle. Higher dimensional
Seifert fibered manifolds were investigated
in \cite{or-wa}. The authors observed that in many cases of interest
Seifert fibered manifolds
 correspond to holomorphic Seifert $\c^*$-bundles, and
started to develop
a general theory of holomorphic Seifert $G$-bundles for any 
complex Lie group $G$.

\begin{defn} Let $X$ be a normal variety (or algebraic space)
over a  field $k$  and $G$ a (reduced)
 algebraic group over $k$. A {\it Seifert $G$-bundle}
over $X$ is a normal variety (or algebraic space) $Y$
together with a morphism $f:Y\to X$ and a 
$G$-action on $Y$
satisfying the following two conditions.
\begin{enumerate}
\item $f$ is affine and $G$-equivariant
(with respect to the trivial action on $X$).
\item For every $x\in X$, the $G$-action on the reduced fiber
$Y_x:=\red f^{-1}(x)$ 
$G\times Y_x\to Y_x$ is $G$-equivariantly isomorphic to
the natural (left) $G$-action on $G/I_x$ for some 
finite subgroup $I_x\subset G$.
\end{enumerate}
\end{defn}

 Let $G\times Y\to Y$ be a proper $G$-action
on a normal variety (or algebraic space). 
The geometric quotient $Y/G$ exists as an algebraic space
(cf. \cite{koll-q, ke-mo}) and 
$f:Y\to Y/G$ is a
Seifert $G$-bundle. As with many quotient problems, even if $Y$
is a quasi projective variety, $X$ need not be quasi projective.
(Many such examples are given in \cite{cox}.)
Thus it is natural to work with algebraic spaces. This is, however,
a purely technical point, and makes no difference for 
our purposes. In the questions motivating this work, $X$ is
always projective.

One can thus view  the theory of 
Seifert $G$-bundles as a special chapter of the
study of algebraic group actions.
The emphasis is, however, quite different.
\cite{or-wa} studied Seifert $S^1$-bundles over $\c\p^n$
in order to construct exotic spheres, extending
and reformulating earlier results of \cite{briesk}.

In a series of papers, Boyer and Galicki developed a method
to construct Einstein metrics on  links of certain
weighted homogeneous singularities \cite{BG01,bgn02,bgn03,BGN03b,bg03,bgk}.
These links are Seifert $S^1$-bundles over the
corresponding projective hypersurfaces in weighted projective
spaces.
 Seifert $S^1$-bundles over smooth rational surfaces
were used to construct Einstein metrics on 
connected sums of $S^2\times S^3$ \cite{ko-s2s3}, but
further results need the study of
Seifert $S^1$-bundles over singular surfaces as well \cite{ko-es5}.

The main aim of this paper is to provide the
necessary foundations on Seifert $\c^*$-bundles
to continue work in this direction, but
Seifert $G_m$-bundles are also interesting from other points of view,
see \cite{dolgachev, pinkham, demazure, fl-za}.

 Ultimately one is mostly interested
in Seifert $S^1$-bundles  $f:M\to X$ where
$M$ is a manifold and $X$  an orbifold, but
the general theory is the same for any normal
variety $X$.

A  complete classification of 
 holomorphic Seifert $\c^*$-bundles
over smooth projective varieties $X$ such that $H_1(X,\z)$ is torsion free
is given in \cite{or-wa}.  In section 1  this is  extend  to
any normal variety   in any characteristic (\ref{classify.thm}).
The local structure of Seifert $G_m$-bundles is studied in section 2.
For applications a key point is the
 smoothness criterion  given in (\ref{good.smooth.crit}).
Section 3 is devoted to investigating the relationship between
holomorphic and topological   Seifert $\c^*$-bundles.
The most interesting Seifert $G_m$-bundles,  those with
a smooth total space, are studied in section 4. Finally Section 5
contains some information about 
the topology of Seifert $\c^*$-bundles. The computation of the
 cohomology groups $H^i(Y,\q)$ is  easy, but the much more interesting
torsion in $H^i(Y,\z)$ remains mostly unexplored.

\section{Classification of Seifert $G_m$-bundles}

\begin{notation} $G_m$ denotes the multiplicative group
$GL(1)$. As a scheme, it is $\spec_k k[t,t^{-1}]$.
The  $M$th roots of unity form the
subgroup scheme 
$$
\mu_M:= \spec_k k[t,t^{-1}]/(t^M-1).
$$
These are all the subgroup schemes of $G_m$. 
$\mu_M$ is nonreduced when $\chr k$ divides $M$.

Every linear representation  $\rho:G_m\to GL(W)$ is completely reducible,
and the same holds for any subgroup $G\subset G_m$
(see, for instance, \cite[I.4.7.3]{groth-sga3-1}).
The set of vectors $\{v:\rho(\lambda)(v)=\lambda^iv\}$
is called the {\it $\lambda^i$-eigenspace}.
We use this terminology also for $\mu_M$-actions.
In this case $i$ is determined modulo $M$.

If a group $G$ acts on a scheme $X$
via $\rho:G\to \aut(X)$, we get an
action on rational functions on $X$ given by
$f\mapsto f\circ \rho(g^{-1})$. (The inverse is needed mostly 
for noncommutative groups only.)

Thus if $G_m$ acts on itself by multiplication,
we get an  induced action on $k[t,t^{-1}]$  where
$\lambda\in G_m(\bar k)$ acts as
$t^i\mapsto \lambda^{-i}t^i$.
Thus $t^i$ spans the $\lambda^{-i}$-eigenspace.

A $G_m$-action on a $k$-algebra $A$ is equivalent to
a $\z$ grading $A=\sum_{i\in \z}A_i$ where 
$A_i$ is  the $\lambda^{-i}$-eigenspace.

The natural $G_m$-action on $G_m/\mu_M$
corresponds to the algebra
$\sum_{i\in M\z}k\cong k[t^M,t^{-M}]$.
\end{notation}

\begin{defn} Let $f:Y\to X$
 be a  Seifert $G_m$-bundle.
 For every $x\in X$, the $G_m$-action on the reduced fiber
$Y_x:=\red f^{-1}(x)$  is isomorphic to
the natural $G_m$-action on $G_m/\mu_{m(x)}$ for some $m(x)$,
called the {\it multiplicity} of the fiber over $x$.
We always assume that  $m(x)=1$ for all $x$ in some dense  open subset of $X$.
That is, we assume that the $G_m$ action is effective.
\end{defn}

\begin{defn} Let $X$ be a normal variety.
A {\it rank 1 reflexive sheaf} is a coherent sheaf $L$
such that $L$ is locally free of rank 1 over a dense open set
and the natural map from $L$ to its double dual
$L\to (L^*)^*$ is an isomorphism.

Rank 1 reflexive sheaves on $X$ form a group, called the {\it class group}
 of $X$,
denoted by $\cl(X)$. The group operation is given by
the double dual of the tensor product $(L\otimes M)^{**}$.
For tensor powers we use the notation
$$
L^{[k]}:=(L^{\otimes k})^{**}.
$$

Let $D$ be a {\it Weil
divisor} on $X$. Then $\o_X(D)$ is a 
 rank 1 reflexive sheaf and we obtain an isomorphism
$$
\cl(X)\cong \weil(X)/(\mbox{linear equivalence}).
$$
For a real number $s$ its {\it round down}
(or  integral part)
is denoted by $\rdown{s}$.
\end{defn}

Our first aim is to classify Seifert $G_m$-bundles over $X$ in terms of
more familiar objects on $X$.

\begin{defn} \label{basic.Gm.exmp}
(cf.\ \cite{pinkham, demazure})
Let $X$ be a normal variety, $L$ a rank 1 reflexive sheaf
on $X$, $D_i$ distinct irreducible  divisors
 and $0< s_i<1$ rational numbers.
Define
$$
\begin{array}{lcl}
S(L, \sum s_iD_i)& :=&\sum_{j\in \z} L^{[j]}
\bigl(\sum_i \rdown{js_i}D_i\bigr)\qtq{and}\\
Y(L, \sum s_iD_i)& :=&\spec_X S(L, \sum s_iD_i).
\end{array}
$$
There is a natural $G_m$-action on $S(L, \sum s_iD_i)$
where $L^{[j]}(\sum_i \rdown{js_i}D_i)$ is the  $\lambda^{j}$ eigensubsheaf.
(See (\ref{sign.choice}) for chosing $\lambda^{j}$
instead of $\lambda^{-j}$.)
This induces a $G_m$-action on $Y(L, \sum s_iD_i)$.

When using the above notation we always assume that the
multiplicative structure of the algebra $S(L, \sum s_iD_i)$ is given by
the tensor product
$$
 L^{[n]}\bigl(\textstyle{\sum_i} \rdown{ns_i}D_i\bigr)\otimes
 L^{[m]}\bigl(\textstyle{\sum_i} \rdown{ms_i}D_i\bigr) \to 
 L^{[n+m]}
\bigl(\textstyle{\sum_i} \rdown{(n+m)s_i}D_i\bigr).
$$
(This map exists since $\rdown{a}+\rdown{b}\leq \rdown{a+b}$ for any $a,b$.)
\end{defn}

\begin{rem}\label{sign.choice}[Sign choices]

Let $f:L\to X$ be a line bundle with zero section
$X\subset L$. One usually identifies
$L$ with its sheaf of sections, temporarily denoted  by
${\mathcal L}$. There is a natural $G_m$-action on
$Y:=L\setminus X$ coming from the standard $GL(1)$ action
on the fibers of $f$. 

The push forward of the structure sheaf is
$f_*\o_L=\sum_{i\geq 0}{\mathcal L}^{-i}$
and ${\mathcal L}^{-i}$ is the $\lambda^{-i}$-eigensubsheaf.
Thus 
 $$
f_*\o_Y=\textstyle{\sum_{i\in\z}{\mathcal L}^i},
$$
and $L$ is naturally identified with ${\mathcal L}$,
which is the $\lambda^1$-eigensubsheaf.
\end{rem}

Special cases of the classification of Seifert $G_m$-bundles
can be found  in the papers 
\cite{dolgachev, pinkham, demazure, fl-za}. The proof of the general
case is essentially the same.

\begin{thm} \label{classify.thm} Let $X$ be a normal variety
(or  algebraic space).
\begin{enumerate}
\item 
Every Seifert $G_m$-bundle $f:Y\to X$ can be uniquely written as 
 $Y\cong Y(L, \sum s_iD_i)$ for some $L$ and $\sum s_iD_i$ as in
(\ref{basic.Gm.exmp}).
\item  $f:Y(L, \sum s_iD_i)\to X$ is a Seifert $G_m$-bundle
iff $L^{[M]}(\sum_i (Ms_i)D_i)$ is locally free for some
$M>0$ and $Ms_i$ is an integer for every $i$.
\item The multiplicity of the Seifert fiber $Y_x$ is the
smallest natural number $M=m(x)$ such that 
 $L^{[M]}(\sum_i (Ms_i)D_i)$ is locally free  at $x$
 and $Ms_i$ is an integer for every $i$ such that $x\in D_i$.
\end{enumerate}
\end{thm}

\begin{rem} (1) One can easily see that the proof also
works for nomal analytic spaces. Thus we obtain that
if $X$ is a normal  and proper variety over $\c$ then
one can naturally identify algebraic and analytic
Seifert $\c^*$-bundles over $X$.

(2)  The result suggests that codimension 1 fixed points
of the $G_m$-action are encoded in the divisorial part
$\sum s_iD_i$ while information about higher codimension fixed point sets is
carried by the sheaf part $L$.
This is, however, probably the wrong interpretation.

From the stacky point of view of (\ref{stacky.rem})
one can identify the pair $(L, \sum s_iD_i)$
with a  single rank one reflexive sheaf on the quotient stack
$Y/G_m$.

(3) Without the local freeness assumption in (\ref{classify.thm}.2),
the finite generation of the algebra $S(L, \sum s_iD_i)$ is a
very subtle question, closely intertwined with the existence of
flips, cf. \cite[Sec.6.1]{kmbook}.

(4) The Theorem  
implies that all Seifert $G_m$-bundles
over $X$ form a  group isomorphic to
$$
\weil(X)_{\q}\times \cl(X)/\langle (D, \o(-D)): D\in \weil(X)\rangle,
$$
but I did not find the group structure useful.
\end{rem}

\begin{defn}\label{chen.class.defn}\cite{or-wa, vist}
Let  $Y\cong Y(L, \sum s_iD_i)\to X$  be a Seifert $G_m$-bundle
over $X$. Then  $L^{[M]}(\sum_i (Ms_i)D_i)$ is locally free 
and so it is an element of $\pic(X)$. We can formally divide
this by $M$ and get an element, called the {\it Chern class}
of $Y/X$
$$
\begin{array}{rrl}
c_1(Y/X)&:=&\tfrac1{M}L^{[M]}(\textstyle{\sum_i} (Ms_i)D_i)\in \pic(X)_{\q}\\
& = & c_1(L)+\textstyle{\sum_i} s_i[D_i].
\end{array}
$$
If we are over $\c$, we also get a topological
Chern class
$$
c_1(Y/X):=\tfrac1{M}c_1\bigl(L^{[M]}(\textstyle{\sum_i} (Ms_i)D_i)\bigr)
\in H^2(X,\q).
$$
There is very little chance of confusion by using the same notation.

Note that although $c_1(Y/X)$ is a rational class,
$M\cdot c_1(Y/X)$ is a well defined class in $H^2(X,\z)$
or in $\pic(X)$. 
More generally, if $j:U\into X$ is an open set and 
$m(U):=\lcm\{m(x):x\in U\}$ then there is a well defined class
$$
m(U)\cdot c_1(Y/X) =
c_1\bigl(j^*\bigl(L^{[m(U)]}(\textstyle{\sum_i} (m(U)s_i)D_i)\bigr)\bigr)
\in H^2(U,\z)
$$
or in $\pic(U)$. 

Similarly, if
 $s(U):=\lcm\{(\mbox{denominator of $s_i$}): D_i\cap U\neq\emptyset\}$ 
then there is a well defined class
$s(U)\cdot c_1(Y/X)\in \cl(U)$.
\end{defn}

\begin{say}[Proof of (\ref{classify.thm}]\label{class.thm.pf}

 Let $f:Y\to X$ be a Seifert $G_m$-bundle.
Since $f:Y\to X$ is affine, $f_*\o_Y$ is a quasicoherent
sheaf with a $G_m$-action. Thus it decomposes as a sum
of $G_m$-eigensubsheaves
$$
f_*\o_Y=\textstyle{\sum_{j\in\z} L_j},
\eqno{(\ref{class.thm.pf}.1)}
$$
where $L_j$ is the  $\lambda^{j}$ eigensubsheaf,
with multiplication maps
$m_{ab}:L_a\otimes L_b\to L_{a+b}$.

Pick any point $x\in X$. 
By assumption $Y_x\cong G_m/\mu_{m(x)}$, thus
$t^{-m(x)}$ on $G_m$ descends to an invertible function $h_x$
on $Y_x$ which is a $G_m$-eigenfunction with eigencharacter
$m(x)$. There is an
  affine neighborhood $x\in U\subset X$ such that 
 $h_x$ lifts to an invertible function $h_U$ 
on $f^{-1}(U)$ which is a $G_m$-eigenfunction with eigencharacter
$m(x)$. This $h_U$ is  a generator of $L_{m(x)}$ on $U$
and $h_U^s$ is a generator of $L_{m(x)}^{\otimes s}$ on $U$.
Thus for  $M=m(X):=\lcm\{m(x):x\in X\}$,  $L_M$ is locally free on $X$
and the multiplication maps
$L_M^{\otimes s}\to L_{sM}$ are isomorphisms for every $s\in \z$.

Taking $M$th power gives a map
$L_1^{\otimes M}\to L_M$, hence
$L_M\cong L_1^{[M]}(\sum d_iD_i)$ where the $D_i$ are
distinct irreducible divisors and $d_i>0$.

For $j>0$, the $j$th power map $L_1^{\otimes j}\to L_j$ 
shows that $L_j\subset (L_j)^{**}= L_1^{[j]}(\sum d_{ij}D_i)$
 for some $d_{ik}$
(where a priori we may have divisors $D_i$ that do not appear 
in the expression for $L_M$ above)
and
the $M$th power map $L_j^{\otimes M}\to L_{jM}\cong L_1^{[jM]}(\sum jd_iD_i)$ 
shows that $\sum d_{ij}D_i\leq \sum \rdown{jd_i/M}D_i$.

This also holds for  $j<0$ as shown by 
 the isomorphisms $L_j\cong L_{j+M}\otimes L_M^{-1}$.

Thus $f_*\o_Y=\sum_{j\in\z} L_j$ is a 
subalgebra of $S(L,\sum (d_i/M)D_i)$ and the two agree
in all degrees divisible by $M$.
Therefore we get that 
$Y(L,\sum (d_i/M)D_i)\to Y$ is 
birational and finite, hence
normalization.
We have assumed that $Y$ is normal, thus
$Y(L,\sum (d_i/M)D_i)\cong Y$.
Since $L_1=L_1^{[1]}(\sum \rdown{d_i/M}D_i)$,
we conclude that $s_i:=d_i/M<1$.
This proves (1).

Conversely, assume that $L^{[M]}(\sum_i (Ms_i)D_i)$ is locally free for some
$M>0$ and $Ms_i$ is an integer for every $i$.
We need to prove that $S(L,\sum s_iD_i)$ is a Seifert $G_m$-bundle.
The question is local, so let us fix a point $x\in X$ and
let $M$ be the smallest  positive number such that
$L^{[M]}(\sum_i (Ms_i)D_i)$ is locally free at $x$ 
 and $Ms_i$ is an integer for every $i$ such that $x\in D_i$.

The fiber of $Y(L,\sum s_iD_i)\to X$ over $x$ is
the spectrum of 
$$
S(L,\textstyle{\sum_i s_iD_i})\otimes k(x)\cong 
\sum_j L^{[j]}(\textstyle{\sum_i (js_i)D_i})\otimes k(x),
$$
where $k(x)$ is the  residue field of $x$.

By (\ref{nilp.prod.lem}), the summands corresponding to
those $j$ such that $L^{[j]}(\sum_i (js_i)D_i)$
is not locally free are nilpotent, and also those
summands where $Ms_i$ is not an integer for every $i$ such that $x\in D_i$.
Thus
\begin{eqnarray*}
Y_x&=&\spec_{k(x)}
\bigl(\sum_j L^{[j]}(\textstyle{\sum_i} (js_i)D_i)\otimes k(x)\bigr)/
(\mbox{nilpotents})\\
&=&
\spec_{k(x)} \sum_{j\in M\z} k(x)\cong G_m/\mu_M.
\end{eqnarray*}
This proves parts (2) and (3).
\qed
\end{say}

\begin{lem} \label{nilp.prod.lem}
Let $L,M$ be rank 1 torsion free sheaves
and assume that there is a surjective map
$h:L\otimes M\onto \o_X$. Then $L,M$ are both locally free.
\end{lem}

Proof. Pick $x\in X$. By assumption there is an
affine nieghborhood $x\in U$ and sections
$\alpha\in H^0(U,L), \beta\in H^0(U,M)$ such that
$h(\alpha\otimes \beta)$ is invertible.

Let $\gamma\in H^0(U,L)$ be arbitrary. Then
$h(\gamma\otimes \beta)=f\cdot h(\alpha\otimes \beta)$
for some $f\in \o_U$, thus
$h((\gamma-f\alpha)\otimes \beta)=0$.
Thus $\gamma-f\alpha$ is zero on the open set where $M$ is locally free,
hence it is zero since $L$ is torsion free. Thus $\alpha$
generates $L|_U$ and so  $L$ is locally free.\qed

\begin{exmp}[Seifert $G_m$-bundles over curves]\label{over.curves.exmp}

Given relatively prime natural numbers $0< b<c$,
let $0<e<c$ be the unique solution of $be\equiv 1\mod c$.

Consider $k[s,s^{-1},u]$ with a $G_m$-action
$s\mapsto \lambda^{-c}s, u\mapsto \lambda^e u$.
The $G_m$-invariants form the polynomial ring $=k[s^eu^c]$.

Set $X:=\spec_k[s^eu^c]$
and
$Y:=\spec_kk[s,s^{-1},u]\to X$. This is a Seifert $G_m$-bundle
with
 $L_c=s^{-1}k[s^eu^c]$
and $L_1=u^bs^{(be-1)/c}k[s^eu^c]$.
Therefore
$$
L_1^{\otimes c}=(u^bs^{(be-1)/c})^ck[s^eu^c]=
(s^eu^c)^bL_c.
$$
Thus, over the line $X$
we get the Seifert $G_m$-bundle
 $Y\cong Y(\o_X, \frac{b}{c}(\mbox{origin}))$. 

The local generating section of $L_1$ is
$h=u^bs^{(be-1)/c}$, thus, as a function on $Y$,
$h$ vanishes along the central fiber with
multiplicity $b$.

Etale locally, this describes Seifert $G_m$-bundles over curves.
More generally, by replacing $k$ with the quotient field
of a divisor $D\subset X$, this construction
describes the behaviour of an arbitrary Seifert $G_m$-bundle
generically (\'etale locally) along divisors.
\end{exmp}

\begin{say}[Maps between Seifert $G_m$-bundles]
\label{seif.maps}

Let $Y=\spec_X\sum_{j\in\z}L_j$ be a Seifert $G_m$-bundle as 
in (\ref{class.thm.pf}.1).
For any natural number $M$,
 $Y_M:=\spec_X\sum_{j\in M\z}L_j$ is also a Seifert $G_m$-bundle,
which can also be viewed as $Y/\mu_M$, the quotient
of $Y$ by the action of $\mu_M\subset G_m$.
Note that $c_1((Y/\mu_M)/X)=Mc_1(Y/X)$.

The case when $M=m(X)$ is especially interesting, as
$\spec_X\sum_{j\in m(X)\z}L_j$ is a $G_m$-bundle.
We denote this by $Y/\mu_X$ with quotient map $\pi:Y\to Y/\mu_X$.

Thus we see that every Seifert $G_m$-bundle can be realized
as a (ramified) cover of a $G_m$-bundle. This observation
was a key step in the classification of \cite{or-wa}.
\end{say}

\begin{say}[Compactification of Seifert $G_m$-bundles]
\label{seif.comp.say}

Seifert $G_m$-bundles can be compactified by adding the
missing zero and infinity sections. Adding the zero section
corresponds to
$\spec_X\sum_{j\leq 0}L_j$ and adding the infinity section
corresponds to
$\spec_X\sum_{j\geq 0}L_j$
in the notation of (\ref{class.thm.pf}.1).
Their union gives a proper morphism
$\bar f:\bar Y\to X$ which is a $\p^1$-bundle over the set where
$m(x)=1$. 
We have the zero section $X_0
\subset \bar Y$ and the infinity section  $X_{\infty}
\subset \bar Y$.
Although $\bar Y$ is almost always singular,
it is a very useful compactification.

For instance, if $X$ is proper, the zero section
$X_0\subset \bar Y$ can be contracted to a point  iff $c_1(Y/X)$ is negative.
The result of the contraction is the affine variety
$\spec_k \sum_{j\leq 0}H^0(X,L_j)$.

Similarly, if $c_1(Y/X)$ is positive, then the infinity section  $X_{\infty}
\subset \bar Y$ can be contracted.

Thus Seifert $G_m$-bundles with $c_1(Y/X)$  negative or positive
can all be viewed as a singularity with  a good $G_m$-action
(minus the singular point itself). This is the point  of view in
\cite{dolgachev, pinkham, demazure}.
\end{say}

\begin{say}[The class group of Seifert bundles]
\label{seif.class.gp}

Let $f:Y=Y(L,\sum\frac{b_i}{a_i}D_i)\to X$ be a Seifert bundle.
Our aim is to compute the class group $\cl(Y)$ of divisors modulo
linear equivalence in terms of $\cl(X),L$ and $\sum\frac{b_i}{a_i}D_i$.

Note first that there is a natural pull back map
$f^*:\cl(X)\to \cl(Y)$ since $f$ is equidimensional.
Set $D^Y_i:=\red f^{-1}(D_i)$.

\begin{prop} \cite[4.22]{fl-za}
 Let  $f:Y=Y(L,\sum\frac{b_i}{a_i}D_i)\to X$ be a Seifert 
$G_m$-bundle.
 Then
$$
\cl(Y)=
\langle f^*\cl(X), D^Y_1,\dots,D^Y_m\rangle/
(f^*[D_i]-a_iD^Y_i, f^*[L]+\sum b_iD^Y_i).
$$
It is somewhat less precise but more transparent to write this as
$$
\cl(Y)\cong 
\langle \cl(X), \tfrac1{a_1}[D_1],\dots,\tfrac1{a_n}[D_n]\rangle/
(c_1(Y/X)).
$$
\end{prop}

Proof. By \cite{fmss}, if a connected solvable group
$G$ acts on a scheme $Y$ then $\cl(Y)$ is generated by
$G$-invariant divisors and the relations are given by
$G$-eigenfunctions.

The $G_m$-invariant divisors on $Y$ are the pull backs of
divisors on $X$, except that $f^*D_i$ becomes $a_iD^Y_i$.

A $G_m$-equivariant rational function on $Y$ can be identified with a
$G$-equivariant rational section of $f_*\o_Y$.
Let $h$ be any nonzero rational section of $L_1$.
Then every  $G_m$-equivariant rational section of $f_*\o_Y$
is of the form $\phi(x)h^k$ for some $k\in \z$,
hence $(h)=0$ gives all other relations.

Choose $h$ such that it has neither zero nor pole along the
divisors $D_i$. As we computed in (\ref{over.curves.exmp}), 
$h$, as a function on $Y$, has a $b_i$-fold zero along
$D^Y_i$. This gives the relation
$f^*[L]+\sum b_iD^Y_i=0$.
\qed  
\medskip

All the rational Chow groups
of arbitrary  Seifert $G$-bundles are computed in \cite[4.4]{vist}.

\begin{say}[Seifert $G_m^n$-bundles] The description of all
Seifert $G_m^n$-bundles for $n\geq 2$ is very similar.

Let $X$ be a normal variety, $L_1,\dots,L_n$  rank 1 reflexive sheaves
on $X$ and  $\Delta_i=\sum_j s_{ij}D_{ij}$  $\q$-divisors
such that $\rdown{\Delta_i}=0$.
Define
$$
\begin{array}{lcl}
S(L_1,\dots, L_n, \Delta_1,\dots, \Delta_n)& :=&
\sum_{j_1,\dots,j_n\in \z} (L_1^{j_1}\otimes\cdots\otimes  L_n^{j_n})^{**} 
\bigl(\rdown{\sum_ij_i\Delta_i}\bigr),\qtq{and}\\
Y(L_1,\dots, L_n, \Delta_1,\dots, \Delta_n)& :=&
\spec_X S(L_1,\dots L_n, \Delta_1,\dots, \Delta_n).
\end{array}
$$
This  is a Seifert $G_m^n$-bundle
iff the set of all multi indices $(j_1,\dots,j_n)$ such that
$$
(L_1^{j_1}\otimes\cdots\otimes  L_n^{j_n})^{**} 
\bigl(\rdown{\textstyle{\sum_i}j_i\Delta_i}\bigr)
$$
is locally free  and $\textstyle{\sum_i}j_i\Delta_i$
is an integral divisor,
is a subgroup of rank $n$ of $\z^n$.

Exactly as in (\ref{classify.thm})
we obtain that 
every Seifert $G_m^n$-bundle $f:Y\to X$ can be  written as 
 $Y\cong Y(L_1,\dots, L_n, \Delta_1,\dots, \Delta_n)$
 for some $L_i$ and $\Delta_i$. 
\end{say}

\section{Local classification of Seifert $G_m$-bundles}

\begin{say}[Quotient description]
Let $f:Y(L,\sum s_iD_i)\to X$ be a Seifert $G_m$-bundle
and $x\in X$ a point. Set $M=m(x)$. By replacing $X$ by a smaller
 open neighborhood of $x$ if necessary, we may assume that
$L_M$ is  free and $L_{jM}\cong L_M^{\otimes j}$ for every $j\in \z$.
Fix an isomorphism $\phi:L_M\cong \o_X$.

$\phi$ defines an algebra structure on
$\sum_{j=0}^{M-1} L_j$ by the rules
$$
\begin{array}{l}
L_a\otimes L_b\stackrel{m_{ab}}{\to} L_{a+b}
\qtq{if $a+b<M$, and}\\
L_a\otimes L_b\stackrel{m_{ab}}{\to} L_{a+b}= L_M\otimes L_{a+b-M}
\stackrel{\phi\otimes 1}{\to}L_{a+b-M} \qtq{if $a+b\geq M$.}
\end{array}
$$
We denote this algebra by $S(L,\sum s_iD_i)/(\phi)$,
the notation suggesting the natural 
 quotient map
$$
\Phi:S(L,\textstyle{\sum s_iD_i})\to S(L,\textstyle{\sum s_iD_i})/(\phi).
$$
Set $Z(L,\sum s_iD_i,\phi):=\spec_XS(L,\sum s_iD_i)/(\phi)$.
It is a Cartier divisor in  $Y(L,\sum s_iD_i)$.

It should be noted that $Z(L,\sum s_iD_i,\phi)$
depends only sligtly on the choice of $\phi$.
Any other choice of $\phi$ can be written as $u\cdot \phi$ where
$u\in \o_X$ is invertible near $x$.
If $\chr k$ does not divide $m(x)$ then $u=v^{m(x)}$ for some
$v$ (at least \'etale locally) and we obtain that
$Z(L,\sum s_iD_i,\phi)\cong Z(L,\sum s_iD_i,u\cdot \phi)$
by the map $L\to L$ which is multiplication by $v$.

There is  a $\mu_M$-action on $S(L,\sum s_iD_i)/(\phi)$
where $L_j$ is 
the $\lambda^j$-eigensubsheaf.

There is also a natural injection
$$
\bar{\Phi}:S(L,\textstyle{\sum s_iD_i})\into  
S(L,\textstyle{\sum s_iD_i})/(\phi)\otimes_k k[t,t^{-1}]
$$
given by $\bar{\Phi}(L_j)=\Phi(L_j)\otimes t^j$.

The $\mu_M$-action on $Z(L,\sum s_iD_i,\phi)$
and the natural  $\mu_M$-action on $G_m$
give a diagonal
$\mu_M$ action on $Z(L,\sum s_iD_i,\phi)\times G_m$
given by
$$
L_i\otimes t^j\mapsto \lambda^{i-j}\cdot L_i\otimes t^j.
$$

\begin{prop}\label{general.local.class.prop}
 The induced map
$\bar \Phi^*:Z(L,\sum s_iD_i,\phi)\times G_m\to Y(L,\sum s_iD_i)$
 is the quotient map by the above $\mu_M$ action. The quotient map
$\bar \Phi^*$  is \'etale if $\chr k$ does not divide $M$.
\end{prop}  

Proof. $L_i\otimes t^j$ is $\mu_M$-invariant iff $i\equiv j\mod M$.
These form exactly the image of $\bar \Phi$.

The $\mu_M$ action on $G_m$ is fixed point free if 
$\chr k$ does not divide $M$, hence so is the 
$\mu_M$ action on $Z(L,\sum s_iD_i,\phi)\times G_m$.
Thus $\bar \Phi^*$  is \'etale.\qed
\end{say}

\begin{prop} \label{gen.smooth.cor}
$Y(L,\sum s_iD_i)$
is smooth along $f^{-1}(x)$ iff
 $Z(L,\sum s_iD_i,\phi)$
is smooth above $x$.
\end{prop}

Proof. If $\chr k$ does not divide $m(x)$ then 
this follows directly from (\ref{general.local.class.prop}).

To see the general case, note that we have computed 
in (\ref{class.thm.pf}) that
the reduced fiber $Y_x$
can be identified as $\spec\sum_{j\in M\z}k(x)$.
Since $S(L,\sum s_iD_i)/(\phi)$ is the sum of terms of degree
les that $M$, we see that $Y_x\subset Y(L,\sum s_iD_i)$ and
the hypersurface $Z(L,\sum s_iD_i,\phi)\subset Y(L,\sum s_iD_i)$
intersect scheme theoretically in one point only.

If $Z(L,\sum s_iD_i,\phi)$ is smooth then so is 
$ Y(L,\sum s_iD_i)$ since  $Z(L,\sum s_iD_i,\phi)$ is  Cartier divisor.
Conversely, in a smooth variety, a Cartier divisor
transversal to  a smooth curve is again smooth.\qed

\begin{say}[Stacky viewpoint]\label{stacky.rem}

Instead of considering the geometric quotient
$Y/G_m$ as a variety, one can also view it is a stack
(see, for instance, \cite{fant-stack} for a short introduction).
 Then we have proved that
the $Z(L,\sum s_iD_i,\phi)/\mu_M$ give local charts
for the quotient stack. This approach will be
especially useful when $Y$ is smooth and $\chr k=0$. In this case
the $Z(L,\sum s_iD_i,\phi)$ are smooth and
so the quotient stack $Y/G_m$ is an orbifold.
This case is studied in detail beginning (\ref{loc.smooth.class}).
\end{say}

While (\ref{gen.smooth.cor}) is a smoothness criterion in principle,
it is   difficult to apply in practice.
Our next aim is to develop a smoothness criterion that is
truly transparent.

It is easier to work from the other direction.
First we recall the local classification of cyclic group actions
near a fixed point, then we see how these give actions on
a $(\mbox{smooth variety})\times G_m$. Then we compute the
corresponding description as $Y(L,\sum s_iD_i)$.

\begin{notation}\label{orbif.notation}
 We consider $\mu_m$-actions on a smooth variety
near a fixed point. We always assume that $\chr k$ does not divide $m$.

Etale locally the action is isomorphic to the induced action
on the tangent space of the fixed point. Thus we need to
classify linear $\mu_m$-actions on $\a^n$ and their quotients.
The action of $\mu_m$  can be diagonalized and described
on $k[z_1,\dots,z_n]$ by $z_i\mapsto \lambda^{a_i}z_i$.
We always assume that $\gcd(a_1,\dots, a_n,m)=1$.
The quotient of $\a^n$ by
 this action of $\mu_m$ is denoted by
$$
\a^n/\mu_m(a_1,\dots,a_n).
$$
(It would be  more consistent to denote this 
by $ \a^n/\mu_m(-a_1,\dots,-a_n)$ instead.
The two actions have the same invariants, so it is probably not
worth while to carry the extra minus sign around.)

 This representation is unique, except that
we are allowed to permute the $a_i$ and 
the $\mu_m$ action can be given by any generator. That is,
if $\gcd(s,m)=1$ then 
$$
\a^n/\mu_m(a_1,\dots,a_n)\cong
\a^n/\mu_m(sa_1,\dots, sa_n).
$$
We write $\a^n_{\mathbf z}$ if we need to indicate that
we are working with an affine space  with coordinates $z_i$.

Given $a_1,\dots,a_n$ and $m$ set
 $$
c_i:=\gcd(a_1,\dots,\widehat{a_i},\dots, a_n,m),\quad
d_i:=a_ic_i/C\qtq{and}C:=\prod c_i.
\eqno{(\ref{orbif.notation}.1)}
$$
Note that the $c_i$ are pairwise relatively prime and $C/c_i$ divides $a_i$.
Observe that $\mu_{c_i}\subset \mu_m$ acts trivially on all but the $i$th
coordinate of $\a^n$,
 so it is a quasi reflection. 
Therefore the quotient
of $\a^n_{\mathbf z}$ by $\mu_C\cong \prod \mu_{c_i}$
is again an affine space $\a^n_{\mathbf x}$ with $x_i=z_i^{c_i}$.
Thus, as a variety,
$$
\a^n_{\mathbf z}/\mu_m(a_1,\dots,a_n)\cong 
\a^n_{\mathbf x}/\mu_{m/C}(d_1,\dots,d_n).
$$
\end{notation}

\begin{say}[The class group of a quotient singularity] 

 Set $X:=\a^n_{\mathbf x}/\mu_M(d_1,\dots,d_n)$
where $\gcd(d_1,\dots, \widehat{d_i},\dots,d_n,M)=1$ for any $i$.
Define furthermore 
$$
D_i: (x_i=0)/\mu_M(d_1,\dots,\widehat{d_i},\dots,d_n)\subset X,
$$
and let
$$
\o_X(j):=\epsilon^j\mbox{- eigenspace of }k[x_1,\dots,x_n],
$$
 as an $\o_X$-module. 
Thus
$$
\o_X(j)=\langle \prod x_i^{v_i}: v_i\geq 0, 
\sum d_iv_i\equiv j\mod M\rangle.
$$
Observe that
$$
\begin{array}{rcl}
\o_X(D_j)&=&\langle \prod x_i^{v_i}: v_j\geq -1,\ v_i\geq 0\ (i\neq j), 
\sum d_iv_i\equiv 0\mod M\rangle\\
&=&x_j^{-1}\langle \prod x_i^{w_i}: w_i\geq 0, 
\sum d_iw_i\equiv d_j\mod M\rangle\cong \o_X(d_j). 
\end{array}
$$

\begin{claim}\label{class.gp.claim}
 The map $j\mapsto \o_X(j)$
gives an  isomorphism $\z/M\cong\cl(X)$. Its inverse is
denoted by $c_1$, called the {\it local Chern class}.
(The ``natural'' Chern class maps $\cl(X)$ to the
group of characters of $\mu_M$. We have identified the
latter with $\z/M$, but it is probably not a good idea
to think of our $c_1$
as a truly  canonical map.)
\end{claim}

Proof. This has been known in various forms for a long time;
for instance the computations of \cite{briesk-s, mumf} both easily
generalize to this case.

 It is probably quickest to note that
the natural $G_m^n$-action on $\a^n$ descends to $X$.
By \cite{fmss}, $\weil(X)$ is generated by $G_m^n$-invariant
divisors, the $D_i$, with $G_m^n$-eigenfunctions giving the
relations.\qed

\end{say}

\begin{say}[Local classification of smooth Seifert bundles]
\label{loc.smooth.class}
From (\ref{general.local.class.prop}) we see that a
Seifert $G_m$-bundle $f:Y\to X$
such that  $Y$ is smooth along $f^{-1}(x)$ is
\'etale locally isomorphic to a Seifert $G_m$-bundle
of the form
$$
f:G_m\times \a^n/\mu_M(r,a_1,\dots,a_n)\to \a^n/\mu_M(a_1,\dots,a_n),
$$
where $r,a_1,\dots,a_n\in \z/M$,  $\gcd(a_1,\dots,a_n,M)=1$.
As before, the Seifert $G_m$-bundle determines the numbers
$r,a_1,\dots,a_n$ modulo $M$, 
except that $(r,a_1,\dots,a_n)$ and  $(sr,sa_1,\dots, sa_n)$
correspond to the same Seifert $G_m$-bundle
for any $s$ which is relatively prime to $M$.

The $\mu_M$-action on the coordinate ring
of $G_m\times \a^n$ is given by
$$
t^s\prod z_i^{u_i}\mapsto \epsilon^{-(rs+\sum a_iu_i)}t^s\prod z_i^{u_i}.
$$
Thus the coordinate ring of $G_m\times \a^n/\mu_M(r,a_1,\dots,a_n)$
is
$$
k[t^s\prod z_i^{u_i}: s,u_i\geq 0, rs+\sum a_iu_i\equiv 0\mod M].
$$
The $G_m$-action is
$$
t^s\prod z_i^{u_i}\mapsto \lambda^{-s}t^s\prod z_i^{u_i}.
$$
In the notation of (\ref{class.thm.pf}.1),
$L_j$ is generated by  monomials of the form $ t^{-j}\prod z_i^{u_i}$, thus
$$
L_j=\langle \prod z_i^{u_i}: u_i\geq 0, \sum a_iu_i\equiv rj\mod M\rangle .
$$

\begin{claim}
The total space $G_m\times \a^n/\mu_M(r,a_1,\dots,a_n)$
is smooth iff either
\begin{enumerate}
\item  $r$ is relatively prime to $M$, or
\item $r=0$ and $\a^n/\mu_M(a_1,\dots,a_n)$ is smooth.
\end{enumerate}
\end{claim}

Proof. If $\chr k\not\vert M$ then the $\mu_M$-action is base point free
if $r$ is relatively prime to $M$, and so the quotient is smooth.
Conversely, if $r\neq 0$, the action of some elment
of $\mu_M$ has a codimension at least 2 fixed point set, so the
quotient can not be smooth by the purity of branch loci.

This settles the problem is characteristic 0.
The positive characteristic cases can all be lifted
to characteristic 0, hence the conditions are necessary.

We still need to show that if  $r$ is relatively prime to $M$
then $G_m\times \a^n/\mu_M(r,a_1,\dots,a_n)$
is smooth. Choose $e_i$ such that $re_i\equiv a_i\mod M$.
Then the coordinate ring of the quotient can be written as
$$
k[x_1t^{-e_1}, \dots,x_nt^{-e_n}, t^M, t^{-M}]
$$
and so $G_m\times \a^n/\mu_M(r,a_1,\dots,a_n)\cong G_m\times \a^n$.\qed
\medskip

We would like to identify $L_1$ as a reflexive sheaf over 
$\a^n/\mu_M(a_1,\dots,a_n)$ and also compute $\sum s_iD_i$.
We start with some remarks on congruences.

\begin{claim}\label{congruence.claim}
 Let $c_i$ and $C$ be as in (\ref{orbif.notation}.1).
\begin{enumerate}
\item One can write
$r$ uniquely as 
$r\equiv lC+\sum a_ib_i\mod M$ where
 $0\leq b_i<c_i$ for every $i$.
\item If $\sum a_iu_i\equiv r\mod M$ then $u_j\equiv b_j\mod c_j$
for every $j$.
\end{enumerate}
\end{claim}

Proof. 
Since $\gcd(a_1,\dots,a_n,M)=1$, we can write
$r\equiv\sum B_ia_i\mod M$. 
$c_i$ divides every $a_j$ except $a_i$, so
$B_i$ modulo $c_i$ is uniquely determined by
$a_iB_i\equiv r\mod c_i$. Let $b_i$ be the remainder of $B_i$ modulo $c_i$.

Assume that $\sum a_iv_i\equiv lC+\sum a_ib_i\mod M$.
Since $c_j|C$ and $c_j|a_i$ for $i\neq j$, we conclude that
$a_jv_j\equiv a_jb_j\mod c_j$. Furthermore, $a_j$ is relatively prime to
$c_j$, hence    $v_j\equiv b_j\mod c_j$. This proves  (2) and the uniqueness
part of (1).
\qed
\medskip

In order to identify $L_1$ as module over
$\a^n/\mu_M(a_1,\dots,a_n)$ we use
(\ref{congruence.claim}.1--2)  to write it as
$$
\begin{array}{lcl}
L_1&=&\langle \prod z_i^{u_i}: u_i\geq 0, \sum a_iu_i\equiv r\mod M \rangle \\
&=&  \bigl(\prod z_i^{b_i}\bigr)
\langle \prod z_i^{u_i-b_i}: u_i\geq b_i, \sum a_iu_i\equiv 
lC+\sum a_ib_i\mod M\rangle \\
&=& \bigl(\prod z_i^{b_i}\bigr)
\langle \prod x_i^{(u_i-b_i)/c_i}: u_i-b_i\geq 0, 
\sum a_ic_i\frac{u_i-b_i}{c_i}\equiv lC\mod M\rangle \\
&=&  \bigl(\prod z_i^{b_i}\bigr)
\langle \prod x_i^{v_i}: v_i\geq 0, 
\sum d_iv_i\equiv l\mod M/C\rangle. \\
\end{array}
$$
Thus $L_1\cong \o_X(l)$.

The divisorial part $\sum s_iD_i$ is computed
from the map $L_1^{\otimes M}\to L_M$.

Since $\gcd(d_1,\dots,\widehat{d_j},\dots, d_n,M/C)=1$ for every $j$, 
the congruence
condition $\sum d_iv_i\equiv l\mod M/C$ has solutions where
$v_j=0$. This implies that the $M$-fold product map
$$
\langle \prod x_i^{v_i}: v_i\geq 0, 
\sum d_iv_i\equiv l\mod M/C\rangle^{\otimes M}\to \o_X
$$
is an isomorphism generically along the $D_j$, hence its kernel and cokernel
are both supported on  a codimension 2
set. Thus we conclude that
$$
L_1^{[M]}=\bigl(\prod z_i^{b_i}\bigr)^ML_M=
\bigl(\prod x_i^{Mb_i/c_i}\bigr)L_M.
$$
Thus the divisorial part of the
representation of our Seifert $G_m$-bundle
is given by $\frac1{M}\sum \frac{Mb_i}{c_i}D_i=\sum \tfrac{b_i}{c_i}D_i$.
We can summarize our computations:

\begin{prop} The Seifert $G_m$-bundle
$$
G_m\times \a^n_{\mathbf z}/\mu_M(r,a_1,\dots,a_n)
\to \a^n/\mu_M(a_1,\dots,a_n)\cong \a^n/\mu_{M/C}(d_1,\dots,d_n)=X
$$
is isomorphic to
$$
Y(\o_X(l),\textstyle{\sum} \tfrac{b_i}{c_i}D_i)\to X,
$$
where  $c_i:= \gcd(a_1,\dots,\widehat{a_i},\dots a_n,M)$
and $l,b_1,\dots,b_n$ are the  unique solutions to
$r\equiv l\prod c_i+\sum a_ib_i\mod M$ satisfying
 $0\leq b_i<c_i$ for every $i$.\qed
\end{prop}

We are ready to formulate a smoothness criterion
that is easy to  use in practice.

\begin{prop}\label{good.smooth.crit}
 Let $X=\a^n/\mu_M(d_1,\dots,d_n)$
be a quotient singularity such that 
$\gcd(d_1,\dots,\widehat{d_i},\dots d_n,M)=1$
for every $i$.
Let $D_i=(x_i=0)/\mu_M(d_1,\dots,\widehat{d_i},\dots,d_n)\subset X$.
Assume that we have  integers $0\leq b_i<c_i$
such that $\gcd(b_i,c_i)=1$ for every $i$.
The following are equivalent.
\begin{enumerate}
\item $Y:=Y(\o_X(l),\sum \tfrac{b_i}{c_i}D_i)$ is smooth.
\item 
\begin{enumerate}
\item  $\gcd(c_i,c_j)=1$ for  every $i\neq j$,
 and 
\item $\o_X(l\prod_i c_i)(\sum_i (\prod_{j\neq i} c_j)b_iD_i)$
 generates $\cl(X)$.
\end{enumerate}
\item
\begin{enumerate}
\item  $\gcd(c_i,c_j)=1$ for  every $i\neq j$, and 
\item $(\prod_i c_i)c_1(Y/X)$ (cf.\ (\ref{chen.class.defn}))
is relatively prime to $M$.
\end{enumerate}
\end{enumerate}
\end{prop}

Proof. In the quotient representation
$Y=G_m\times \a^n_{\mathbf z}/\mu_m(r,a_1,\dots,a_n)$
the $c_i$ are pairwise relatively prime and 
$d_i$ is relatively prime to $c_i$. Furthermore, 
$Y$ is smooth iff
$r=(\prod c_i)l+\sum_i (\prod_{j\neq i} c_j)b_id_i$ is 
relatively prime to $m=M\prod c_i$.

Pick a prime $p$ dividing $c_k$. Then $p\not\vert b_kd_k$
and $p$ divides every summand of 
$(\prod c_i)l+\sum_i (\prod_{j\neq i} c_j)b_id_i$ save
$ (\prod_{j\neq k} c_j)b_kd_k$. Thus $r$ being relatively
prime to $ m=M\prod c_i$ is equivalent to
being relatively
prime to $M$.

The conditions (b) also imply that  $\gcd(c_i,d_i)=1$
for every $i$.
\qed
\end{say}

\section{Topological Seifert $\c^*$-bundles}

\begin{defn} \cite{or-wa} \label{top.seif.Gbund.defn}
Let $X$  
 be a nomal analytic space
 and $G$ a real Lie group.

A {\it topological Seifert $G$-bundle} over 
$X$ is a topological space
$Y$ together with a 
 $G$-action and a
continuous  map $f:Y\to X$
such that  
$X$ has an open covering $X=\cup_iU_i$
such that for every $i$ 
  the preimage $f:f^{-1}(U_i)\to U_i$ is
fiber preserving $G$-equivariantly homeomorphic to
a ``standard local Seifert $G$-bundle''
$$
f_i:(G\times V_i)/G_i\to U_i.
$$
Here $V_i$ is a normal analytic space, $G_i$ a finite group
of biholomorphisms  such that $V_i/G_i\cong U_i$
and
 the $G_i$-action on $G\times V_i$ is the diagonal action
given 
on $G$ by a homomorphism
$\phi_i:G_i\to G$ composed with  the right action of $G$ on itself.
 The left $G$-action of $G$ on itself
gives the $G$ action on $Y$.

In order to avoid 
nontrivial orbifold structures on $Y$, we always assume
that the $G_i$-action on $G\times V_i$ is fixed point free outside
a codimension 2 set.

Every analytic  Seifert $G$-bundle corresponds to a
topological Seifert $G$-bundle by ignoring the analytic structure.
\end{defn}

\begin{rem} From the topological point of view it is rather unnatural
to require that the $V_i$ be complex spaces and it may be interesting to
develop a purely topological theory without any analytic assumptions.
Our main interest, however, lies in understanding the analytic and
topological properties of holomorphic Seifert $\c^*$-bundles.

Aside from some smoothness questions, the above definition
is equivalent to the one in \cite{or-wa}.
\end{rem}

\begin{lem}\label{local.top=anal} Let
 $f_i:Y_i\to X$  be two holomorphic Seifert $\c^*$-bundles
which are $\c^*$-equivariantly homeomorphic.
Then there is an open covering $X=\cup U_j$ such that
$f_1^{-1}(U_j)\to U_j$ and $f_2^{-1}(U_j)\to U_j$
are $\c^*$-equivariantly biholomorphic for every $j$.
\end{lem}

Proof. For sufficiently divisible $M$, consider the quotient 
Seifert $\c^*$-bundles
$Y_i/\mu_M\to X$ defined in (\ref{seif.maps}).

Both are
$\c^*$-bundles, hence there is a covering $X=\cup U_j$
such that the $Y_i/\mu_M\to X$ are trivial over every $U_j$.

The quotient maps $Y_i\to Y_i/\mu_M$ are finite degree ramified
coverings which are $\c^*$-equivariantly homeomorphic to each other.

The biholomorphisms  $Y_1/\mu_M|_{U_j}\cong Y_2/\mu_M|_{U_j}$
differ from the global homeomorphism $Y_1/\mu_M\sim Y_2/\mu_M$
by a translation by a continuous function $U_j\to \c^*$. Thus by
suitably changing the $\c^*$-equivariant homeomorphism
over $U_j$ we get a commutative diagram
$$
\begin{array}{ccc}
Y_1|_{U_j} & \sim &  Y_2|_{U_j}\\
\downarrow && \downarrow\\
Y_1/\mu_M|_{U_j}& \cong & Y_2/\mu_M|_{U_j},
\end{array}
$$
where $\sim$ is a homeomorphism, $\cong$ a 
biholomorphism and the vertical arrows are
finite degree ramified
coverings of analytic spaces. Thus 
$Y_1|_{U_j}$ is  biholomorphic to $Y_2|_{U_j}$
by the uniqueness of analytic structures on finite coverings.\qed
\medskip

The following rather standard  result provides an alternate
approach to the description of Seifert $\c^*$-bundles.
It also provides an efficient way to compare
holomorphic and topological Seifert $\c^*$-bundles.

\begin{defn} 
 Let $X$ be a normal analytic space
and  $G$  a commutative complex Lie group.
Holomorphic sections of $G\times X\to X$ form a sheaf, denoted by $G^{an}_X$.
Its
cohomology groups are denoted by $H^i(X,G^{an})$.
Similarly,  continuous sections of $G\times X\to X$ form a sheaf,
 denoted by $G^{top}_X$. Its
cohomology groups are denoted by $H^i(X,G^{top})$.
There are natural ``forgetful'' maps  $H^i(X,G^{an})\to H^i(X,G^{top})$.

Note that $(\c^*)^{an}_X$ is usually denoted by $\o_X^*$.
\end{defn}

\begin{prop} \label{seif.obstr}
 Let $X$ be a normal analytic space
and  $G$  a commutative complex (resp.\ real) 
Lie group. Let $X=\cup U_i$ be an open cover 
and
assume that over each $U_i$ we have an analytic (resp.\ topological)
 Seifert $G$-bundle $Y_i\to U_i$
and  there are $G$-equivariant biholomorphisms
(resp.\  homeomorphisms) 
$\phi_{ij}:Y_j|_{U_{ij}}\cong Y_i|_{U_{ij}}$.
\begin{enumerate}
\item There is an obstruction element in the torsion subgroup
 $H^2_{tors}(X,G^{an})$
(resp. $H^2_{tors}(X,G^{top})$)
such that there is a global Seifert $G$-bundle $Y\to X$ compatible with these
local structures iff the obstruction element is zero.
\item The set of all such global Seifert $G$-bundles,  up to 
$G$-equivariant biholomorphisms (resp.\ homeomorphisms),
 is either empty or forms a principal homogeneous
space under $H^1(X,G^{an})$
(resp.\ $H^1(X,G^{top})$).
\item Both of the constructions commute with the
forgetful maps $H^i(X,G^{an})\to H^i(X,G^{top})$.
\end{enumerate}
\end{prop}

Proof. 
The analytic and topological proofs are exactly the same.
We write $H^i(X,G)$ to indicate that we can work in
either settings.

The isomorphisms $\phi_{ij}$ can be changed
to  $\alpha_{ij}\phi_{ij}$ for any
$\alpha_{ij}\in H^0(U_{ij},G)$.  These patchings define a 
global Seifert $G$-bundle iff
$$
\alpha_{ik}\phi_{ik}=\alpha_{ij}\phi_{ij}\alpha_{jk}\phi_{jk}
\qtq{for every $i,j,k$.}
$$
This is equivalent to
$$
\alpha_{ij}\alpha_{jk}\alpha_{ki}=(\phi_{ij}\phi_{jk}\phi_{ki})^{-1}
\qtq{for every $i,j,k$.}
\eqno{(\ref{seif.obstr}.4)}
$$
The products $(\phi_{ij}\phi_{jk}\phi_{ki})^{-1}\in H^0(U_{ijk}, G)$
satisfy the cocycle condition, and they define
an element of $H^2(X,G)$, called the obstruction.
One can find
$\{\alpha_{ij}\}$ satisfying
(\ref{seif.obstr}.4) iff the obstruction is zero.

Replacing the $Y_j$ by $Y_j/\mu_M$ changes the
isomorphisms over $U_i\cap U_j$ to $\phi_{ij}^M$, hence the
obstruction corresponding to the Seifert $\c^*$-bundles 
$Y_j/\mu_M$ is the $M$th power of the original obstruction.

If $M$ is sufficiently divisible, the quotients
 $Y_j/\mu_M$ are all $\c^*$-bundles, and these can always be
globalized to the trivial $\c^*$-bundle. Thus the obstruction is
always torsion.

Two choices  $\{\alpha_{ij}\}$ and $\{\alpha'_{ij}\}$
give isomorphic Seifert bundles iff there are
isomorphisms $\delta_i:Y_i\cong Y_i$
(viewed as elements of $H^0(U_i,G)$)
such that
$$
\alpha'_{ij}\alpha^{-1}_{ij}=\delta_i\delta^{-1}_j|_{U_{ij}}.
$$
Thus  $\{\alpha'_{ij}\alpha^{-1}_{ij}\}$
corresponds to a class in $H^1(X,G)$.\qed

\begin{exmp} Consider $\p^1\times \p^1$. Blow up the 4 points all
 of whose coordinates are in $\{0,\infty\}$ and contract the
resulting 4 rational curves of selfintersection $-2$.
We get a surface $S$ with 4 ordinary double points $p_i$.
The local class group of the points is $\z/2$.

By looking at all curves on $\p^1\times \p^1$ we easily see that
if $C\subset S$ is any curve and $[C]_i$ is the corresponding element
in the  local class group of $p_i$ then
$\sum_i [C]_i=0\in \z/2$. In particular we obtain that any
Seifert $\c^*$ bundle which is trivial on $S\setminus \{p_1\}$
is trivial.

Here $H^3(S,\z)\cong \z/2$ and there is a nontrivial obtsruction.
\end{exmp}

\begin{lem} \label{reduce.to.ord.coh} Notation as above.
\begin{enumerate}
\item $H^i(X,(\c^*)^{top})\cong H^{i+1}(X,\z)$ for $i\geq 1$.
\item If  $H^i(X,\o_X)=0$ for $i=1,2$  then 
$H^1(X,(\c^*)^{an})\cong H^2(X,\z)$  and 
there is an injection $H^2(X,(\c^*)^{an})\into H^3(X,\z)$.
\end{enumerate}
\end{lem}

Proof. Let $\c^{top}_X$ denote the sheaf of continuous complex valued
 functions. This sheaf is soft and so it has no higher cohomologies
(cf.\ \cite[II.9]{bredon}).
Thus the 
 long exact cohomology sequence of the 
exponential sequence
$$
0\to \z\stackrel{2\pi i}{\to} \c^{top}_X\stackrel{exp}{\to} (\c^*)^{top}_X
\to 0
$$ 
proves (1).

The second part follows similarly from the holomorphic
exponential sequence
$$
0\to \z\stackrel{2\pi i}{\to} \o_X\stackrel{exp}{\to} \o_X^*\to 0.\qed
$$ 

\begin{rem}
 Let $X$ be a normal analytic space with rational
singularities such that $H^i(X,\o_X)=0$ for $i=1,2$.

 Let $f:X'\to X$ be any resolution of singularities.
It is not hard to see (cf. \cite{ar-mu}) that 
$\ker \left[H^3_{tors}(X,\z)\stackrel{f^*}{\to} H^3_{tors}(X',\z)\right]$
is independent of the choice of $X'$.

One can show that the obstructions in (\ref{seif.obstr}.1)
are in this smaller group. We do not use this in the sequel.
\end{rem}

\begin{thm}  Let $X$ be a normal analytic space
such that $H^i(X,\o_X)=0$ for $i=1,2$. Then 
the forgetful map
$$
\left\{
\begin{array}{c}
\mbox{holomorphic}\\ 
\mbox{Seifert $\c^*$-bundles}
\end{array}
\right\}
\longrightarrow
 \left\{
\begin{array}{c}
\mbox{topological}\\ 
\mbox{Seifert $\c^*$-bundles}
\end{array}
\right\}
$$
is an isomorphism.
\end{thm}

Proof. First we consider injectivity. By (\ref{local.top=anal})
if $f_i:Y_i\to X$  are two holomorphic Seifert $\c^*$-bundles
which are homeomorphic then they are
 locally biholomorphic to each other.

Thus, by (\ref{seif.obstr})
 they differ by an element of $H^1(X,(\c^*)^{an})$,
which is isomorphic to $H^2(X,\z)$ by (\ref{reduce.to.ord.coh}). 
Hence they are also different as topological Seifert $\c^*$-bundles,
again 
by  (\ref{seif.obstr}) and (\ref{reduce.to.ord.coh}). 

In order to prove that every topological Seifert $\c^*$-bundle
has a holomorphic structure, consider the quotient 
Seifert $\c^*$-bundle
$Y/\mu_X$ defined in (\ref{seif.maps}).
(Although these were considered in the algebraic case, the construction can
be also performed topologically.)

The quotient $Y/\mu_X$ is a topological $\c^*$-bundle, and
these are classified by $H^2(X,\z)$. By our
assumptions and (\ref{reduce.to.ord.coh}), every
toplogical $\c^*$-bundle has a unique analytic $\c^*$-bundle
structure.

The quotient map $Y\to Y/\mu_X$ is a finite degree ramified
covering map, and it ramifies along the preimages of
the branch divisors of $V_i\to U_i$ which are analytic
subvarieties. Thus $Y$ has a unique
analytic structure which makes the quotient map $Y\to Y/\mu_X$
holomorphic. \qed

\section{Orbifolds}

\begin{defn}[Orbifolds]\label{orbif.defn}
  An {\it orbifold}  is 
 a normal variety $X$ covered by
\'etale charts $\cup U_i\onto X$ where each $U_i$
is
written as  a quotient of a smooth variety by a finite group.
That is,  for each $U_i$ there
is a smooth variety $V_i$ and a finite (reduced) group
$G_i$ acting on $V_i$ such that $U_i$ is
isomorphic to the quotient space $V_i/G_i$.
The quotient maps are denoted by $\phi_i:V_i\to U_i$.

We  denote an orbifold by  $X^{orb}$.

The compatibility condition between the charts
that one needs to assume is that 
the coordinate projections from the normalization of the fiber products
$$
\overline{V_i\times_XV_j}\to V_i
$$
are all \'etale.

To avoid complications, we only consider the case
when the orders $|G_i|$ are not divisible by the characteristic.

Since the transition maps between the charts are \'etale,
the cotangent bundles $\Omega^1_{V_i}$
glue together to the orbifold cotangent bundle
$\Omega^1_{X^{orb}}$. 

Although we write $X=\cup U_i$, we are always allowed to replace the
covering $\{U_i\}$
by a suitable finer covering $\{U_{i_j}\}$.

We allow the case  when there are 
codimension 1 fixed point sets. 
 Then the quotient map $\phi_i:V_i\to U_i$ has 
branch divisors $D_{ij}\subset U_i$
and ramification divisors $R_{ij}\subset V_i$.
Let $m_{ij}$ denote the ramification index
over $D_{ij}$. Locally near a general point of
$R_{ij}$ the map $\phi_i$ looks like
$$
\phi_i:(x_1,x_2,\dots,x_n)\mapsto 
(z_1=x_1^{m_{ij}},z_2=x_2,\dots,z_n=x_n).
$$
The compatibility condition between the charts
implies that  there are
global divisors $D_j\subset X$ and ramification
indices $m_j$ such that $D_{ij}=U_i\cap D_j$
and $m_{ij}=m_j$ (after suitable re-indexing).

It will be convenient to codify these data
by a single $\q$-divisor, called the {\it branch divisor}
of the orbifold,
$$
\Delta:=\sum (1-\tfrac1{m_j})D_j.
$$

The orbifold is
uniquely determined by the pair $(X,\Delta)$.
Indeed, given a point $v\in V_i$, set $u=\phi_i(v)\in U_i$.
Then \'etale locally the quotient map
$\phi_i: (v\in V_i)\to (u\in U_i)$ is uniquely determined
by the following condition:
\begin{enumerate}
\item $\phi$ is unramified over $U_i\setminus(\sing U_i\cup\bigcup_j D_j)$.
\item The ramification index of $\phi$ over $D_j$ divides $m_j$. 
\item $\phi_i: (v\in V_i)\to (u\in U_i)$ is the maximal covering satisfying
the above 2 conditions.
\end{enumerate}

Slightly inaccurately, we  sometimes  identify the orbifold
with the pair $(X,\Delta)$.

An orbifold $(X,\Delta)$ 
is called {\it locally cyclic} if  all the local charts
$U_i\subset X$ are quotients of the $V_i$ by cyclic groups $G_i$.
Thus, \'etale locally, it can be given by charts
$$
\a^n_{\mathbf z}/\mu_M(a_1,\dots,a_n),
$$
where
$\gcd(a_1,\dots,a_n,M)=1$.

Let $f:Y\to X$ be a Seifert $G_m$-bundle and assume that $\chr k=0$.
By (\ref{gen.smooth.cor}), this defines an orbifold structure
on $X$, denoted by $(X,\Delta)$.
Note that the divisors appearing in
$\Delta=\sum (1-\tfrac1{m_j})D_j$
are the same as the divisors
appearing in the representation
$Y=Y(L,\sum \frac{b_j}{c_j}D_j)$. 
Furthermore, $c_j=m_j$ but the orbifold does not carry
any information about the $b_j$ or about $L$.

We write  a Seifert $G_m$-bundle as
$$
f:Y\to (X,\Delta)
$$
if we want to emphasize the orbifold structure
on the base.

From the point of view of \cite{or-wa, ko-es5} it is very natural
to choose first $(X,\Delta)$ and then consider the
effect of various choices of $b_i$ and $L$ 
later.
\end{defn}

\begin{prop}\label{Omega.of.Y}
 Let $f:Y\to (X,\Delta)$ be a  Seifert $G_m$-bundle
over a field of characteristic 0, $Y$ smooth.
There is an exact sequence
$$
0\to f^*\Omega^1_{X^{orb}}\to \Omega^1_Y\to \o_Y\to 0.
$$
\end{prop}

Proof. For any  orbifold chart 
$U=V/\mu_M$ of $(X,\Delta)$ we get an
\'etale chart 
$G_m\times V/\mu_M$ on $Y$.
Choose  coordinates
$(t, x_1,\dots,x_n)$, then 
$$
 \Omega^1_{G_m\times V}\cong \tfrac{dt}{t}\o_{G_m\times V}+f^*\Omega^1_V.
$$
This decomposition is invariant under the $\mu_M$-action, but
the representation of $G_m\times V$
 as a direct product is not unique, as we can replace
$t$ by $\phi(x)t$ where $\phi\in \o_V$ is any invertible function.
Since
$$
\frac{d(\phi(x)t)}{\phi(x)t}=\frac{dt}{t}+\sum_i \frac1{\phi(x)}
\frac{\partial \phi(x)}{\partial x_i}dx_i,
$$
we conclude that $\tfrac{dt}{t}$ gives a well defined
global generator of $\Omega^1_Y/f^*\Omega^1_{X^{orb}}$.\qed

\begin{cor}\label{K.of.Y} \cite{fl-za}
 Let $f:Y\to (X,\Delta)$ be a  Seifert $G_m$-bundle
over a field of characteristic 0. Then
$$
K_Y=f^*K_X+\sum (a_i-1)D^Y_i.
$$
\end{cor}

Proof. Linear equivalence of divisors is not affected
by throwing away a subset of codimension 2, thus
we may assume that $Y$ and $X$ are smooth.
This case follows directly from (\ref{Omega.of.Y})
by noting that $f^*K_{X^{orb}}=f^*K_X+\sum (a_i-1)D^Y_i$.
\qed
\medskip

This may seem a little unexpected since for
the total space $g:Z\to X$ of a line bundle $L$ the canonical
bundle formula is $K_Z=g^*(K_X+L)$. We are, however, looking 
only at $Y:=Z\setminus(\mbox{zero section})$, and this
exactly kills the $g^*L$ part of the formula.
\end{say}

\begin{say}  Let $f:Y\to (X,\Delta)$ be a  Seifert $G_m$-bundle
over a field of characteristic 0 such that $c_1(Y/X)$ is ample.
As noted in (\ref{seif.comp.say}), this implies that the section at
infinity $X_{\infty}\subset \bar Y$ is contractible
and we obtain a singularity $0\in W$ with a good $G_m$-action.

From (\ref{K.of.Y}) and (\ref{seif.class.gp}) we conclude that
$K_W$ is $\q$-Cartier iff $K_X+\Delta$ is a rational multiple of
$c_1(Y/X)$.

Furthermore, if $(X,\Delta)$ is an orbifold,
then  applying the generalized adjunction formula
(as developed in \cite[Sec.16]{k-etal})  to $X_{\infty}\subset \bar Y$
we conclude that $0\in W$ is  log terminal (cf.\ \cite[2.34]{kmbook})
iff
$-(K_X+\Delta)$ is ample.
\end{say}

\section{Topology  of Seifert $\c^*$-bundles}

In this section the base field is $\c$.

Our aim 
is to obtain information about 
 the integral cohomology groups and  the fundamental group
of a Seifert bundle
$f:Y\to (X,\Delta)$ in terms of
$(X,\Delta)$ and the Chern class  of $Y\to X$.

The cohomology groups $H^i(Y,\z)$ are computed by a 
Leray  spectral sequence whose $E_2$ term is
$$
E^{i,j}_2=H^i(X, R^jf_*\z_{Y})\Rightarrow H^{i+j}(Y,\z).
$$
Every fiber  of $f$ is $\c^*$, so 
$R^2f_*\z_Y=0$ and
the only
interesting
higher direct image is $R^1f_*\z_{Y}$.
Our first task is to compute this sheaf and its
cohomology groups. 

Next we consider the simplest edge homomorphism in the spectral
sequence
$$
\delta: H^0(X, R^1f_*\z_{Y})\to H^{2}(X,\z),
$$
and identify it with  the Chern class  $m(X)c_1(Y/X)$.

In some cases of interest, these data completely determine the
cohomology groups, and even the topology,  of $Y$.
Some of these instances are discussed in
\cite{or-wa, ko-s2s3}.

\begin{prop}\label{R^1.prop} Let $f:Y\to X$ be a Seifert $\c^*$-bundle.
\begin{enumerate}
\item There is a natural isomorphism
$\tau:R^1f_*\q_Y\cong \q_X$.
\item There is a natural injection
$\tau:R^1f_*\z_Y\into \z_X$ which is an isomorphism over points where $m(x)=1$.
\item If $U\subset X$ is connected then
$$
\tau(H^0(U,R^1f_*\z_Y))= m(U)\cdot  H^0(U,\z)\cong m(U)\cdot \z,
$$
where $m(U)$ is the $\lcm$ of the multiplicities of
all fibers over $U$.
\end{enumerate}
\end{prop}

Proof. Pick $x\in X$ and a contractible neighborhood
$x\in V\subset X$. Then $f^{-1}(V)$ retracts to
$S^1\subset f^{-1}(x)$ and (together with the orientation)
this gives  a distinguished generator
$\rho\in H^1(f^{-1}(V), \z)$. This in turn determines
a  cohomology class $\frac1{m(x)}\rho\in H^1(f^{-1}(V), \q)$.
These normalized cohomology classes  are compatible with each other
and give a global section of $R^1f_*\q_Y$. Thus
$R^1f_*\q_Y=\q_X$ and we also obtain the
injection $\tau:R^1f_*\z_Y\into \z_X$ as in (2). 

If $U\subset X$ is connected, a section $b\in \z\cong H^0(U,\z_U)$ is in
$\tau(R^1f_*\z_Y)$ iff $m(x)$ divides $b$ for every $x\in U$.
This is exactly (3).\qed

\begin{cor}\label{chernclass=edgemap.cor}
 The quotient map $\pi:Y\to Y/\mu_X$ (defined in (\ref{seif.maps}))
induces an isomorphism
$$
H^0(X, R^1f_*\z_{Y})\cong \pi^* H^0(X, R^1(f/\mu_X)_*\z_{Y/\mu_X}).
$$
Under this isomorphism, the edge homomorphism
$$
\delta: H^0(X, R^1f_*\q_{Y})\to H^{2}(X,\q)
$$
is identified with the Chern class  $c_1(Y/X)$.
Thus the image of
$$
\delta: H^0(X, R^1f_*\z_{Y})\to H^{2}(X,\z)
$$
is generated by   $c_1(Y/\mu_X)=m(X)c_1(Y/X)$.
\end{cor}

Proof.  If $m(X)$ is the lcm of the multiplicities of all fibers
then 
$$
\tau(\pi^* H^0(X, R^1(f/\mu_X)_*\z_{Y/\mu_X}))=m(X)\cdot H^0(X,\z)
$$
and so it agrees with $\tau(H^0(X, R^1f_*\z_{Y}))$.

The map $\pi$ induces a map between the  spectral sequences
$$
H^i(X, R^jf_*\z_{Y/\mu_X})\Rightarrow H^{i+j}(Y/\mu_X,\z)
\qtq{and}
H^i(X, R^jf_*\z_{Y})\Rightarrow H^{i+j}(Y,\z),
$$
thus $\delta$ is identified with the corresponding
edge homomorphism
$$
\delta: H^0(X, R^1(f/\mu_X)_*\z_{Y/\mu_X})\to H^{2}(X,\z).
$$
For a $\c^*$-bundle this is exactly the 
Chern class of $Y/\mu_X$.\qed

\begin{cor}  The quotient map $\pi:Y\to Y/\mu_X$
induces an isomorphism on  rational cohomologies.
Thus $Y$ is a rational homology sphere iff
$$
H^*(X,\q)\cong \q[c_1(Y/X)]/(c_1(Y/X)^{n+1})\qtq{where 
$n=\dim X$.}
$$
\end{cor}

Proof.  We have seen that $\pi$ induces an isomorphism
betwen the spectral sequences with $\q$-coefficients,
giving the first claim. 
A circle bundle $M\to X$ is an integral (resp.\ rational)
homology sphere iff $c_1(M/X)$ generates $H^*(X,\z)$
(resp.\ $H^*(X,\q)$).\qed

\begin{note} The only
simply connected smooth projective variety 
such that $H^*(X,\z)$ is generated by a degree 2 element that I know is $\p^n$.
All Seifert bundles $Y\to \p^n$ homeomorphic to
$S^{2n+1}$ were described in \cite{or-wa}.

On the other hand, there are many other Seifert bundles
over singular varieties $Y\to X$ such that
 $Y$ is homeomorphic to
$S^{2n+1}$. These give  interesting examples of
Einstein  metrics on $S^{2n+1}$, see \cite{bgk}.
\end{note}

 We see  that (\ref{R^1.prop}) describes
the sheaf $R^1f_*\z_Y$ completely in terms of $(X,\Delta)$,
but  it is not always easy to compute its cohomologies 
based on this description. There are, however, some cases where
this is  quite straightforward.

\begin{prop}\label{smooth.R^1.prop}
 Let  $f:Y\to (X,\sum (1-\frac1{m_i})D_i)$ 
be a Seifert bundle and assume that $X$ is smooth.
Then there is an exact sequence
$$
0\to R^1f_*\z_Y\stackrel{\tau}{\to} \z_X \to \sum_i \z_{D_i}/m_i\to 0.
$$
\end{prop}

Proof. Note that $m_i$ and $m_j$ are relatively prime if 
$D_i\cap D_j\neq \emptyset$. 
It is now clear that the kernel of 
$\z_X \to \sum_i \z_{D_i}/m_i$ has the same sections as
described in (\ref{R^1.prop}.3).\qed
\medskip

We get a slightly more complicated description
if $X$ has only isolated singularities.

\begin{prop}\label{isol.sings.R^1.prop}
 Let  $f:Y\to (X,\sum (1-\frac1{m_i})D_i)$ 
be a Seifert bundle. Assume that $X$ has only
isolated singularities  of type $\a^n/\mu_{n_j}$ at $P_j\in X$ as in
(\ref{orbif.notation}).
Then there is an exact sequence
$$
0\to R^1f_*\z_Y\stackrel{\tau}{\to} \z_X \to Q\to 0.
$$
and in turn $Q$ sits in another exact sequence
$$
0\to \sum_j \z_{P_j}/n_j\to  Q \to \sum_i \z_{D_i}/m_i\to 0.\hfill\qed
$$
\end{prop}

\begin{defn} Let $(X,\sum_i (1-\frac1{m_i})D_i)$ be an orbifold
and  $X^0\subset X$  the smooth locus of $X$.
The {\it orbifold fundamental group} $\pi^{orb}_1(X,\Delta)$
is the fundamental group of $X^0\setminus\supp\Delta$
modulo the relations: if $\gamma$ is any 
 small loop around $D_i$ at a general point
 then $\gamma^{m_i}=1$.

Note that $\pi^{orb}_1(X,0)$ is usually different from
$\pi_1(X)$.

The abelianization of $\pi^{orb}_1(X^0,\Delta)$
is  denoted by  $H^{orb}_1(X^0,\Delta)$, called the
{\it abelian orbifold fundamental group.}
(I do not know how to define analogous orbifold homology groups
$H_i^{orb}$ for $i\geq 2$.)
\end{defn}

The following is a straightforward generalization of the
computation of the fundamental group of 3--dimensional Seifert 
bundles, cf.\ \cite{seif}.

\begin{prop}\label{fund.gr.sequence}
 Let $f:Y\to (X,\Delta)$ be a Seifert bundle
and let $X^0\subset X$ denote the smooth locus. 
Assume that $Y$ is smooth. There is
an exact sequence
$$
\pi_1(S^1)\to \pi_1(Y)\to \pi^{orb}_1(X,\Delta)\to 1.
$$
\end{prop}

Proof. Removing a real codimension 4 set $f^{-1}(X\setminus X^0)$
from $Y$ does not change the fundamental group, thus we may assume
that $X$ is smooth.

We have a surjection 
$$
f_*:\pi_1(f^{-1}(X\setminus\supp\Delta))\onto 
\pi_1(X\setminus\supp\Delta).
$$
 Let  $\gamma_i$ be a 
 small loop around $D_i$ at a general point. 
Then $\gamma_i$ lifts to 
a   loop around $f^{-1}(D_i)$
and the lifting of  $\gamma_i^{m_i}$ contracts in $Y$.
Thus $f_*$ descends
 to a  surjection
$$
f_*:\pi_1(Y)\onto 
\pi^{orb}_1(X,\Delta).
$$
In order to compute the kernel, we can pass to the
universal orbifold cover of $(X,\Delta)$ and
so we assume that $\pi^{orb}_1(X,\Delta)=1$.

Let $\gamma\in Y$ be a closed loop. We can perturb it
and assume that it is contained in $Y\setminus f^{-1}(\supp\Delta)$.
The image $f(\gamma)\in X\setminus\supp\Delta$
need not be contractible, but it is in the normal subgroup
generated by the above  $\gamma_i^{m_i}$. The latter are all
images of loops that are contractible in $Y$, thus by composing
$\gamma$ with these we may assume that $f(\gamma)$ is contractible in
$X\setminus\supp\Delta$. Since $Y\setminus f^{-1}(\supp\Delta)
\to X\setminus\supp\Delta$ is a circle bundle, 
$\gamma$ is thus homotopic to a multiple of the fiber.
\qed

\begin{say}\label{orbif.H_1.say}
The determination of $\pi^{orb}_1(X,\Delta)$
seems rather tricky in general, but  its abelianization
 is
fully computable.
The case when $H_1(X^0,\z)=0$ is especially easy to state.
\end{say}

\begin{prop}\cite[4.6]{or-wa}\label{OW.h1.lem}
 Let $X$ be a complex manifold such that $H_1(X,\z)=0$ and let 
$D_1,\dots,D_n\subset X$ be smooth divisors intersecting transversally.
\begin{enumerate}
\item $H^{orb}_1(X,\sum(1-\frac1{c_i})D_i)$
is given by generators $g_1,\dots,g_n$ and relations
\begin{enumerate}
\item $c_ig_i=0$ for $i=1,\dots,n$, and
\item $\sum g_i([D_i]\cap \eta)=0$
for every $\eta\in H_2(X,\z)$.
\end{enumerate}
\item For any line bundle $L$, $H_1(Y(L,\sum\frac{b_i}{c_i}D_i),\z)$
is given by generators $k,g_1,\dots,g_n$ and relations
\begin{enumerate}
\item $c_ig_i+b_ik=0$ for $i=1,\dots,n$, and
\item $k(c_1(L)\cap \eta)-\sum g_i([D_i]\cap \eta)=0$
for every $\eta\in H_2(X,\z)$.\qed
\end{enumerate}
\end{enumerate}
\end{prop}

\begin{prop} \label{c_1.gives.Sb.for.simpconn}
Yet $(X,\Delta)$ be an orbifold and assume that
$H^{orb}_1(X,\Delta)=0$. Then a Seifert $\c^*$-bundle
$f:Y\to (X,\Delta)$ is uniquely determined 
by its Chern class
$c_1(Y/X)\in H^2(X,\q)$.
\end{prop}

Proof. Let $Y^j\to (X,\Delta)$ 
be two Seifert $\c^*$-bundles
constructed from the invariants  $(\sum (b_i^j/c_i)D_i, L^j)$
 as in (\ref{classify.thm}).
 The equality of their
 Chern classes means that
$$
c_1(L^1)+\sum \tfrac{b_i^1}{c_i}[D_i]=
c_1(L^2)+\sum \tfrac{b_i^2}{c_i}[D_i]\in H^2(X,\q).
$$
Let $X^0\subset X$ denote the smooth locus and set $M:=\lcm\{c_i\}$.
Since $H_1(X^0,\z)=0$, there is no torsion in $H^2(X^0,\z)$,
and so after restricting to $X^0$ we get an equality in integral
cohomology
$$
M\cdot (c_1(L^1)-c_1(L^2))+\sum (b_i^1-b_i^2)\tfrac{M}{c_i}[D_i]=0
\in H^2(X^0,\z).\eqno{(\ref{c_1.gives.Sb.for.simpconn}.1)}
$$
$H_1(X,\z)=0$ also implies that $H^1(X,\o_X)=0$, hence the Picard group
of $X^0$ injects into $H^2(X^0,\z)$, and so
(\ref{c_1.gives.Sb.for.simpconn}.1) is also a linear equivalence.

Set $\beta_i=b_i^1-b_i^2$ if $b_i^1-b_i^2\geq 0$ and
$\beta_i=c_i+b_i^1-b_i^2$ if $b_i^1-b_i^2< 0$. Then
(\ref{c_1.gives.Sb.for.simpconn}.1)  can be rearranged to
$$
\sum \beta_i\tfrac{M}{c_i}[D_i]=M\cdot c_1(N),
$$
where $N$ is some line bundle on $X^0$ and $0\leq \beta_i<c_i$.
This corresponds to an $M$-sheeted covering of
$X^0$ which ramifies along $D_i$ with ramification index  
$=c_i/\gcd(\beta_i,c_i)$.
This covering represents a nontrivial element of 
$H^{orb}_1(X,\Delta)$, unless $c_i|\beta_i$ for every $i$. This is only
possible if $\beta_i=0$ for every $i$. Thus $b_i^1=b_i^2$ for every $i$
and hence also $L^1\cong L^2$.\qed

\begin{ack} I thank Ch.\ Boyer, K.\ Galicki, S.\ Keel, T.\ Pantev,
 A.\ Vistoli and M.\ Zaidenberg
for many useful conversations, e-mails and references.
Research
was partially supported by the NSF under grant number
DMS-0200883. 
\end{ack}

\bibliography{refs}

\vskip1cm

\noindent Princeton University, Princeton NJ 08544-1000

\begin{verbatim}kollar@math.princeton.edu\end{verbatim}

\end{document}